\input amstex
\magnification=1095
\input amsppt.sty

\hsize=14truecm \vsize=22.6truecm
\hsize=16.2truecm
\vsize=22truecm

\TagsOnRight
\def\t{\tilde}
\def\W{\bold W}
\def\bv{\bf Vir}
\def\lra{\longrightarrow}
\def\vn{\varepsilon}
\def\ot{\otimes}
\def\b{\overline}

\def\U{\bold U}
\leftheadtext{N. H. Hu} \rightheadtext{$q$-Witt Algebras,
$q$-Virasoro algebra, $q$-Lie Algebras and Representations}

\topmatter

\

\

\

\title\nofrills{$q$-Witt Algebras, $q$-Virasoro algebra, $q$-Lie Algebras, $q$-Holomorph
Structure and Representations}
\endtitle
\author{\bf Naihong Hu}\endauthor

\footnote" "{The work was supported in part by the National Natural Science
Foundation of China (No: 19731004) and the Alexander von Humboldt Foundation.}
\affil
Department of Mathematics, East China Normal University\\
Shanghai 200062, China\\
E-mail: \ nhhu\@math.ecnu.edu.cn
\endaffil

\abstract For $q$ generic or $q=\vn$ a primitive $l$-th root of $1$,
$q$-Witt algebras are described by means of $q$-divided power
algebras. The structure of the universal $q$-central extension of
the $q$-Witt algebra, the $q$-Virasoro algebra $\text{\bv}^q$, is
also determined. $q$-Lie algebras are investigated and the $q$-PBW
theorem for the universal enveloping algebras of $q$-Lie algebras is
proved. A realization of a class of representations of the $q$-Witt
algebras is given. Based on it, the $q$-holomorph structure for the
$q$-Witt algebras is constructed, which interprets the realization
in the context of representation theory.


\endabstract

\endtopmatter

\baselineskip11pt
\noindent
{\bf 1 \ Introduction and Notation}
\vskip0.3cm

\noindent
Much of the recent work in quantum group theory has focused on the
investigation of the structure and representations of the
quantum groups $R_q[G]$ for the semisimple Lie groups $G$ or the quantized
universal enveloping algebras $U_q(\frak g)$ for the classical semisimple (or
affine) Lie algebras $\frak g$. As well-known, especially for the primitive
root of unity case, this leads to an interesting representation theory
mirroring in a remarkable fashion the prime characteristic case for both
algebraic groups and classical Lie algebras (cf.\cite{5},
etc.). In the discussion, $q$-integers and $q$-binomial numbers play
important roles. From the viewpoint of Lie theory, however, the attention to
quantum groups has been paid mainly on those Lie algebras associated to
(generalized) Cartan matrices. The interest here is to wonder how certain
``quantum phenomenon'' might be related to those Lie algebras of Cartan type.
We begin with such an attempt for a simple situation of Witt type.

\vskip0.15cm
As mentioned above, $q$-integers is of basic importance in quantum phenomenon.
The Jackson differential operator as skew derivation was used by Kassel
(\cite{9}) and other authors. It actually makes a close connection between the
Jackson differential operator and $q$-integers. As known in the modular Lie
algebra theory, the divided power structure over a commutative algebra plays a
crucial role in the structure and representation theory of Lie algebras of
Cartan type (\cite{14}, \cite{11--13}). So we hope that there exists a
similar divided power structure reflecting the structure nature of the quantum
version of the Witt algebras. Starting from Jackson operators, we
introduce a $q$-divided power structure in terms of $q$-binomial numbers
such that it is compatible with the characterization of skew derivations.
Thereby, a natural development can be proceeded by following the approach to
modular Lie algebras of Cartan type. As in the modular case,
such a divided power algebra structure not only describes intrinsically the
structure of $q$-Witt algebras, but plays a nice role in a realization of a
class of representations as well. A further insight into the realization
leads to the construction of the $q$-holomorph structure, which interprets
the intrinsic characteristic of the realization in representation theory, so
that more representations (not necessarily graded) can be obtained by it.

\vskip0.15cm
The paper is organized as follows. Section 2 is about a description of
quantum version of the Witt algebras. A weighted $q$-Jacobi identity under a $q$-Lie
product also holds for them. $q$-Holomorph structure for $q$-Witt algebras,
which still obeys the weighted $q$-Jacobi identity, is constructed. These
examples allow us to extract a general algebraic object, $q$-Lie algebras,
in section 3. The $q$-PBW theorem for their universal enveloping algebras is
also valid. As an application of the weighted $q$-Jacobi identity, the structure of
the universal $q$-central
extension of the $q$-Witt algebra, the $q$-Virasoro algebra $\text{\bv}^q$, is
also determined. In section 4, the behavior of the universal enveloping algebra
of the $\vn$-Witt algebra $\W^\vn(1,\bold 1)$ over its central subalgebra is
understood and, a conjecture of de Concini and Procesi is connected with our
case. Section 5 concerns with a realization of a class of representations for
$\W^\vn(1,\bold 1)$, which is established over a tensor vector space
$\frak A^\vn(1,\bold 1)\ot_{\Cal K} V$. By virtue of it, the irreducible
modules are easily described. An analogous realization for the $q$-Witt
algebra $\W^q(1)$ can be transferred automatically to a tensor vector space
$\frak A^q(1)\ot_{\Cal K} V$. A more general nature for such a realization
is further revealed in section 6 by taking advantage of the $q$-holomorph
structure $\Cal H^\vn(1)$ of the $\vn$-Witt algebra $\W^\vn(1,\bold 1)$ with
its associated quantum base $l$-space $_1\Cal K_\vn^l$.

\vskip0.2cm
\noindent
{\it Notation.} \
Let $\Cal K$ be an arbitrary field, $\text{\it char}(\Cal K)\ne 2,3$,
and $q\in \Cal K$, $q\not=0,1$.
For $n, r\in\Bbb Z^+$ ($n\ge r$), set
$$
(n)_q=\frac{1-q^n}{1-q},\quad
(n)_q!=(1)_q(2)_q\cdots (n)_q,\quad
\binom{n}r_q=\frac{(n)_q!}{(n-r)_q!(r)_q!}.
$$
Then $(-n)_q:=-q^{-n}(n)_q$.
We make the convention: $(0)_q!=1$, $\binom{n}0_q=1$
and $\binom{n}r_q=0$ for $n<r$.
It is well known that the $q$-binomial number
$\binom{n}r_q$ is a Gaussian
polynomial in the variable $q$. Therefore, $\binom{n}r_q$ is well defined
for all $q\in\Cal K$.
Notice also that if $q=\vn$ is a primitive root of
unity of order $l$ ($l\in \Bbb N$), then
$(kl)_\vn=0$ and $\binom{l}i_\vn=0$, for $k\in \Bbb Z, k\ne 0$ and $0<i<l$.

\vskip0.5cm
\noindent
{\bf 2 \ Quantum Version of the Witt Algebras}
\vskip0.3cm

\noindent
{\it $q$-Witt algebra $\W^q$.} \
Recall that the Witt algebra $\text{\bf W}=\text{Der}$ $(\Bbb C[x,x^{-1}])$
consists of derivations of the algebra $\Bbb C[x,x^{-1}]$ with Lie product
$$
\split
[x^{i+1}\partial,x^{j+1}\partial]&=x^{i+1}\partial(x^{j+1})\partial-
x^{j+1}\partial(x^{i+1})\partial\\
&=(j-i)x^{i+j+1}\partial,
\endsplit\tag{1}
$$
where $\partial=d/dx$, $i, j\in\Bbb Z$. It is well-known that
$\W(1)=\{x^{i+1}\partial \mid i\in\Bbb Z, i\ge -1\}$ is the Witt subalgebra
of Cartan type in $\W$.

Let $\Cal K$ be an arbitrary field, $\text{\it char}(\Cal K)\ne 2,3$,
and $q\in \Cal K$, $q\ne 0,1$ be generic.
Following \cite{9}, we define Jackson's $q$-differential operator $\partial_q$
over $\Cal K[x,x^{-1}]$ by
$$
\partial_q(\Cal P)=\frac {\Cal P(qx)-\Cal P(x)}{qx-x},\quad \forall
\Cal P\in\Cal K[x,x^{-1}].\eqno(2)
$$

Let $\tau_q$ denote an algebra automorphism of $\Cal K[x,x^{-1}]$ defined by
$\tau_q(x)=qx$. The $q$-differential operator $\partial_q$, which is a linear
mapping over $\Cal K$, is a {\it $\tau_q$-derivation} or {\it skew
derivation} (cf. \cite{5}, \cite{9}), namely for all $\Cal P, \Cal Q \in
\Cal K[x,x^{-1}]$, we have
$$
\partial_q(\Cal {PQ})=
\partial_q(\Cal P)\Cal Q
+\tau_q(\Cal P)\partial_q(\Cal Q).\eqno(3)
$$
Let $\text{Der}_q(\Cal K[x,x^{-1}])$ denote the set of all
$\tau_q$-derivations over $\Cal K[x,x^{-1}]$, and
let $e_n:=x^{n+1}\partial_q$, then we have

\proclaim{Lemma 2.1} $(\text{\rm i})$ \ $\text{Der}_q(\Cal K[x,x^{-1}])$ is a free
$\Cal K[x,x^{-1}]$-module of rank $1$ with $\partial_q$ as its
base

\item{}\qquad \;over $\Cal K[x,x^{-1}]$, where $\partial_q(x^n)=(n)_qx^{n-1}$.

$(\text{\rm ii})$ \ \;$\text{Der}_q(\Cal K[x,x^{-1}])$ is a vector space over $\Cal K$
with a basis $\{ e_n\mid n\in\Bbb Z \}$.

$(\text{\rm iii})$ \ If we define a $q$-bracket product $\{,\}_q$ on
$\text{Der}_q(\Cal K[x,x^{-1}])$ by
$$
\{e_i,e_j\}_q(x^n):=q^{i+1}e_i(e_j(x^n))-q^{j+1}e_j(e_i(x^n)),\eqno(4)
$$

\item{}\qquad
then we have
$$\split
\{e_i,e_j\}_q&=x^{i+1}\partial_q(x^{j+1})\partial_q
-x^{j+1}\partial_q(x^{i+1})\partial_q\\
&=[(j+1)_q-(i+1)_q]\,e_{i+j},\quad(i\in\Bbb Z)
\endsplit\tag{5}
$$

\item{}
\qquad and then the $q$-bracket product $\{,\}_q$ is bilinear over
$\Cal K$ and satisfies the antisymmetry and

\item{}
\qquad the weighted $q$-Jacobi identities:
$$
\gather
\{e_i,e_j\}_q=-\{e_j,e_i\}_q,\tag {6}\\
(2)_{q^i}\{e_i,\{e_j,e_k\}_q\}_q{+}
(2)_{q^j}\{e_j,\{e_k,e_i\}_q\}_q{+}
(2)_{q^k}\{e_k,\{e_i,e_j\}_q\}_q=0.\tag{7}
\endgather
$$
\endproclaim
\demo{Proof} (i) \ The first claim of (i) is clear since the ring
$\Cal K[x,x^{-1}]$ is commutative and for any $\Cal P\in \Cal K[x,x^{-1}]$,
$\Cal P\partial_q$ also enjoys the relation (3). According to
(2), we have
$$\partial_q(x^n)=\frac{(qx)^n-x^n}{qx-x}=(n)_qx^{n-1}.$$

(ii) is obvious due to (i).

(iii) \ By (i), we have
$$\split
&\{e_i,e_j\}_q(x^n)\\
&\quad=(n)_q[q^{i+1}(n+j)_q-q^{j+1}(n+i)_q]\,x^{n+i+j}\\
&\quad=[(j+1)_q-(i+1)_q]\,e_{i+j}(x^n)\\
&\quad=[x^{i+1}\partial_q(x^{j+1})\partial_q
       -x^{j+1}\partial_q(x^{i+1})\partial_q](x^n).
\endsplit
$$
The bilinearity of the $q$-bracket $\{,\}_q$ follows from the linearity
of $\partial_q$, as does the usual Witt algebra. The antisymmetry
(6) is clear according to its definition (5). The weighted $q$-Jacobi identity
(7) is obtained from the following identity and a cyclic permutation
of $(i,j,k)$.
$$\split
(&2)_{q^i}\{e_i,\{e_j,e_k\}_q\}_q\\
&=(1+q^i)[(k+1)_q-(j+1)_q][(j+k+1)_q-(i+1)_q]\,e_{i+j+k}\\
&=\frac{q^2}{1{-}q}[(i{+}k)_q{+}(2i{+}k)_q{+}(2j{+}k)_q{+}(i{+}2j{+}k)_q\\
&\quad -(i{+}j)_q{-}(2i{+}j)_q{-}(j{+}2k)_q{-}(i{+}j{+}2k)_q]\,e_{i{+}j{+}k}.
\endsplit
$$
Thus we complete the proof.\hfill\qed
\enddemo

\noindent
{\it Remark.} \ The definition (5) of $q$-bracket product $\{,\}_q$ is
quite similar to (1), whose bilinearity is naturally ensured.
Let $\W^q=(\text{Der}_q(\Cal K[x,x^{{-}1}]),\{,\}_q)$.
Obviously, when $\Cal K=\Bbb C$, the ordinary Witt algebra
$\bold W=\text{Der}(\Bbb C[x,x^{-1}])$
is the limit of $\bold W^q$ when $q$ tends to $1$.
$\W^q$ is referred to as the $q$-Witt algebra, which
is a $q$-Lie algebra with $q$-Lie product $\{,\}_q$ in the sense of
Definitions 3.1.

\vskip0.3cm
\noindent
{\it $q$-Witt algebra $\W^q(1)$.} \
Let $q\in \Cal K$, $q\ne 0, 1$ and be generic. Associated to the Witt
subalgebra $\W(1)$ of Cartan type in $\W$, an intrinsic description of
$q$-Witt algebra $W^q(1)$ can be given by introducing a so-called
{\it $q$-divided power structure} in terms of $q$-binomial numbers. An
advantage is that the ``skew'' nature of skew derivations is well-revealed
over it (compare (10) with (12)).

\proclaim{Definition 2.2}
An algebra $\frak A^q(1)$ with the generators $\{x^{(a)}\mid a\in\Bbb Z^+\}$
over $\Cal K$ is called {\it a $q$-divided power algebra}, if it has the
multiplication
$$
x^{(a)}x^{(b)}=\binom{a+b}a_q x^{(a+b)},\quad\forall
\ a, b\in \Bbb Z^+.\eqno(8)
$$
\endproclaim

Similar to the modular case (cf. \cite{14}), we define a special
$\tau_q$-derivative of the $q$-divided power algebra $\frak A^q(1)$ by:
$$
\partial_q (x^{(a)})=x^{(a-1)}, \eqno(9)
$$
where $\tau_q$ is an automorphism of $\frak A^q(1)$ defined by
$\tau_q(x^{(a)}):=q^ax^{(a)}$.
Denote $\text{Der}_q(\frak A^q(1))$ by the set of all the special
$\tau_q$-derivatives over $\frak A^q(1)$,
set $e_{(n)}:=x^{(n+1)}\partial_q$.
Similar to Lemma 2.1, we have
\proclaim{Corollary 2.3} $(\text{\rm i})$ \ $\partial_q$ is a special $\tau_q$-derivative
of the $q$-divided power algebra $\frak A^q(1)$, i.e.,
$$
\partial_q(x^{(a)}x^{(b)})=\partial_q(x^{(a)})x^{(b)}
+\tau_q(x^{(a)})\partial_q(x^{(b)}).\eqno(10)
$$

$(\text{\rm ii})$ \ \;$\text{Der}_q(\frak A^q(1))$ is a vector space over $\Cal K$
with a basis $\{e_{(i)}\mid i\ge -1\}$.

$(\text{\rm iii})$ \ The $q$-Lie product $\{,\}_q$ on $\text{Der}_q(\frak A^q(1))$
defined as
$$
\split
\{e_{(i)},e_{(j)}\}_q:
&=x^{(i+1)}\partial_q(x^{(j+1)})\partial_q
 -x^{(j+1)}\partial_q(x^{(i+1)})\partial_q \\
&=\left[\binom{i{+}j{+}1}{i{+}1}_q
 -\binom{i{+}j{+}1}{j{+}1}_q\right]e_{(i{+}j)},\quad (\forall i, j\ge -1)
\endsplit
$$

\item{}\qquad
is bilinear and satisfies the antisymmetry and the weighted $q$-Jacobi identities.
\endproclaim
\demo{Proof} \ (10) follows from
the operation property of $q$-integers and $q$-binomial numbers:
$$
\gather
(a)_q+q^a(b)_q=(a+b)_q,\tag{11}\\
\binom{a+b-1}{a-1}_q+q^a\binom{a+b-1}a_q=\binom{a+b}a_q.\tag{12}
\endgather
$$

(ii) \& (iii) are clear from Lemma 2.1.
\hfill\qed
\enddemo

Denote $\W^q(1):=(\text{Der}_q(\frak A^q(1)),\{,\}_q)$.
It is a $q$-Lie algebra (in the sense of Definitions 3.1).

\vskip0.2cm
\noindent
{\it Remark.} \
Let $q$ be a non-root of unity. Let $e_n:=(n+1)_q!e_{(n)}$, then
it follows from Corollary 2.3 (iii) that $\W^q(1)$ is a $q$-Witt subalgebra
of $\W^q$. In particular, for the $q$-divided power algebra $\frak A^q(1)$
in this case, there is a natural interpretation, i.e.
$x^{(a)}:=\frac{1}{(a)_q!}x^a,\quad \forall a\in \Bbb Z^+$.
Since $\partial_q (x^a)=(a)_q x^{a-1}$, the defining relation $(9)$
is naturally satisfied. In the following, the divided power structure is
crucial to the case of primitive $l$-th root of unity.

\vskip0.3cm
\noindent
{\it $\vn$-Witt algebra $\W^\vn(1,\bold 1)$.} \
Let $q=\vn\in\Cal K$ be a primitive $l$-th root of $1$,
for $l\in \Bbb N$, we consider {\it an $\vn$-divided power subalgebra
$\frak A^\vn(1,\bold 1)$} of $\frak A^\vn(1)$ over $\Cal K$ as follows.
$$
\frak A^\vn(1,\bold 1):=\langle x^{(a)}\mid 0\le a<l, a\in \Bbb Z\rangle,
$$
where $x^{(a)}x^{(b)}=\binom{a+b}{a}_\vn x^{(a+b)}$ and $(x^{(a)})^{l}=0$.

Similarly, we may define the special $\tau_\vn$-derivatives
$e_{(n)}=x^{(n+1)}\partial_\vn$ ($-1\le n\le l-2$) over the $\vn$-divided
power algebra $\frak A^\vn(1,\bold 1)$ and the $\vn$-Lie product $\{,\}_\vn$
of them as in Corollary 2.3. It is clear that
$$
\W^\vn(1,\bold 1):=
\text{Der}_\vn(\frak A^\vn(1,\bold 1))=
\langle e_{(i)}\mid -1\le i\le l-2\rangle
$$
is closed under the $\vn$-Lie product $\{,\}_\vn$, namely,
$$
\{e_{(i)},e_{(j)}\}_\vn=\cases
\left[\binom{i+j+1}{i+1}_\vn-\binom{i+j+1}{j+1}_\vn\right]e_{(i+j)},
& {-}1\le i{+}j\le l{-}2,\\
0, & \text{\it otherwise}.
\endcases\eqno(13)
$$

\noindent
{\it Remark.} \ In the case of primitive $l$-th root of unity, the
commutator relation (13) is quite analogous to that of the
modular Witt algebra $\W(1,\bold 1)$ (cf. \cite{14}).
$\W^\vn(1,\bold 1)$ under $\{,\}_\vn$ satisfies the
antisymmetry and weighted $\vn$-Jacobi identities as in Lemma 2.1.
$\W^\vn(1,\bold 1)$ is called an $\vn$-Witt algebra
over $\Cal K$. In particular, $\W^\vn(1,\bold 1)$ is an $l$-dimensional
$\vn$-Witt subalgebra of $\W^\vn(1)$ (in the sense of Definitions 3.1), which
is similar to the modular case (cf. \cite{14}).

\vskip0.3cm
Recall the notion of quantum space introduced by Manin (\cite{10}).
Let $\bold q=(q_{ij})$ be an $n\times n$ matrix with the entries in $\Cal K^*$
and $q_{ii}=1$,  $q_{ij}^{-1}=q_{ji} \ (i\ne j)$.
Let $\Cal K\{L_1,\cdots,L_n\}$ denote the free algebra generated by
$L_1,\cdots,L_n$ and let $\Cal I_q$ be the ideal of $\Cal K\{L_1,\cdots,L_n\}$
generated by the elements $L_iL_j-q_{ij}L_jL_i$. A {\it quantum $n$-space}
is defined as its quotient-algebra
(with the same symbols $L_i$ writting for the images under the natural homomorphism)
$$\Cal K_q^n:=\Cal K_q[L_1,\cdots,L_n]=\Cal K\{L_1,\cdots,L_n\}/\Cal I_q.$$
The quantum $n$-space $\Cal K_q^n$ is too big for our purpose, we need only
its quantum base space
$_1\Cal K_q^n:=\langle L_i\mid 1\le i\le n
\rangle$ satisfiying the relations $L_iL_j=q_{ij}L_jL_i$.

\vskip0.15cm
The following two special quantum spaces will be needed later.

\vskip0.2cm
\noindent
{\it Example 2.4 $($Quantum $\infty$-space $\Cal K_q^\infty\,)$} \ Take a
$\bold q$-matrix $(q^{j{-}i})_{i,j\in\Bbb Z^+}$, where
$q\in \Cal K^*$ is a non-root of unity.
Denote  by $\Cal K_q^\infty=
\Cal K_q[L_0,L_1,\cdots]$ the $\infty$-quantum space with the relations
$L_iL_j=q^{j-i}L_jL_i \quad (i,j\in \Bbb Z^+)$.

\vskip0.2cm
\noindent
{\it Example 2.5 $($Quantum $l$-space $\Cal K_\vn^l\,)$} \ Let
$q=\vn$ be a primitive
$l$-th root of unity. Take an $\bold \vn$-matrix
$(\vn^{j-i})_{i,j\in \Bbb Z_{(l)}}$,
where $\Bbb Z_{(l)}=\Bbb Z/(l)$.
The quantum $\l$-space $\Cal K_\vn^l$ is denoted by
$\Cal K_\vn^l=\Cal K_\vn[L_0,L_1,\cdots,L_{l -1}]$
with the relations
$L_iL_j=\vn^{j-i}L_jL_i \quad(0\le i,j\le l-1)$.

\vskip0.3cm
\noindent
{\it $q$-Holomorphs.} \
For the $q$-Witt
algebra $\W^q(1)$ and the $\vn$-Witt algebra $\W^\vn(1,\bold 1)$, we
construct their {\it $q$-holomorph structure}, which is a kind of
$q$-extensions of $\W^q(1)$ and $\W^\vn(1,\bold 1)$ through the associated
quantum base spaces respectively (see Definitions 3.1 (iv)). Such objects look like
the ``holomorph" structure in Lie algebra theory, which play a special
role in representations of $\W^q(1)$ and $\W^\vn(1,\bold 1)$ (see
section 6).

Note that the quantum base space $_1\Cal K_q^\infty$ or
$_1\Cal K_\vn^l$ is $q$- or $\vn$-abelian respectively, namely,
$$\{_1\Cal K_q^\infty,\,_1\Cal K_q^\infty\}_q=0,
\quad \{_1\Cal K_\vn^l,\,_1\Cal K_\vn^l\}_\vn=0.$$

Denote
$$\split
\Cal H^q:&=\,_1\Cal K_q^\infty\bigoplus\W^q(1),\\
\quad\Cal H^\vn(1):&=\,_1\Cal K_\vn^l\bigoplus\W^\vn(1,\bold 1).
\endsplit\tag{14}
$$

Define
$$
\split
\{e_{(i)},L_j\}_q&=q^{i+1}\binom{i+j}{i+1}_qL_{i+j}, \quad (
i\ge -1, j\ge 0, i,j \in \Bbb Z),\\
\{e_{(i)},L_j\}_\vn&=\vn^{i+1}\binom{i+j}{i+1}_\vn L_{i+j}, \quad (
 -1\le i\le \l-2, 0\le j\le l-1).
\endsplit\tag{15}
$$
\proclaim{Theorem 2.6} \
The direct sum $\Cal H^\vn(1)$ $($resp. $\Cal H^q)$ forms a new $\vn$-Lie
algebra $($resp. $q$-Lie algebra$)$ under the $\vn$-Lie product $\{,\}_\vn$ $($resp.
the $q$-Lie product $\{,\}_q)$, called the $q$-holomorph structure of
$\Cal H^\vn(1)$ $($resp. $\Cal H^q)$. In particular,
$L_0$ is a $\vn$-central $($resp. $q$-central\,$)$ element in $\Cal H^\vn(1)$
$($resp. $\Cal H^q)$, i.e. $\{L_0,\Cal H^\vn(1)\}_\vn=0$
$($resp. $\{L_0,\Cal H^q\}_q=0)$.
\endproclaim
\demo{Proof} \
For the $\vn$-Witt algebra $\W^\vn(1,\bold 1)$ and its associated quantum
base $l$-space $_1\Cal K_\vn^l$,
we need only to show that
$\Cal H^\vn(1)$ is still an $\vn$-Lie algebra in the sense of Definitions 3.1.

The antisymmetry of $\{,\}_\vn$ for $\Cal H^\vn(1)$ is obvious.

Note that $\{_1\Cal K_\vn^l,\,_1\Cal K_\vn^l\}_\vn=0$.
It suffices to verify that the weighted $\vn$-Jacobi identity holds for
the elements $e_{(i)}, e_{(j)}, L_k$.

From the following identities
$$
\split
&(2)_{\vn^i}\{e_{(i)},\{e_{(j)},L_k\}_\vn\}_\vn \\
&\quad=\frac{\vn^{i+j+2}(i+j+k)_\vn!}{(i+1)_\vn!(j+1)_\vn!(k-1)_\vn!}
(1+\vn^i)(j+k)_\vn L_{i+j+k},\\
&(2)_{\vn^j}\{e_{(j)},\{L_k,e_{(i)}\}_\vn\}_\vn\\
&\quad=-\frac{\vn^{i+j+2}(i+j+k)_\vn!}{(i+1)_\vn!(j+1)_\vn!(k-1)_\vn!}
(1+\vn^j)(i+k)_\vn L_{i+j+k},\\
&(2)_{\vn^k}\{L_k,\{e_{(i)},e_{(j)}\}_\vn\}_\vn\\
&\quad=(1{+}\vn^k)\left[\binom{i{+}j{+}1}{j{+}1}_\vn{-}
\binom{i{+}j{+}1}{i{+}1}_\vn\right]\vn^{i{+}j{+}1}\binom{i{+}j{+}k}{i{+}j{+}1}_\vn
L_{i{+}j{+}k}\\
&\quad=\frac{\vn^{i+j+1}(i+j+k)_\vn!}{(i+1)_\vn!(j+1)_\vn!(k-1)_\vn!}
\left[(i+1)_\vn-(j+1)_\vn\right](1+\vn^k)L_{i+j+k},
\endsplit
$$
it is clear that the weighted $\vn$-Jacobi identity holds for them.

The proof for $\Cal H^q$ is similar.
\hfill\qed
\enddemo

\proclaim{Corollary 2.7} The direct sum $\langle L_0\rangle\bigoplus \W^q(1)$
$($resp.$\langle L_0\rangle \bigoplus \W^\vn(1,\bold 1)\,)$ is a split
$1$-dimensional $q$-central $($resp. $\vn$-central\,$)$ extension of $\W^q$
$($resp. $\W^\vn(1,\bold 1)\,)$.   \hfill\qed
\endproclaim

\noindent
{\it Remark.} \
This case suggests it might be interesting to consider the universal
$q$-central extension of the $q$-Witt algebra $\W^q$. But indeed, in that
case, the extension $\langle L_0 \rangle\bigoplus \W^q$ is not split any
more, just like that relationship between the usual Witt algebra $\W$ and the
Virasoro algebra $\text{\bv}$ (see section 3).

\vskip0.5cm
\noindent
{\bf 3 \ $q$-Lie Algebras, $q$-$\text{PBW}$ Theorem and $q$-Virasoro Algebra}

\vskip0.3cm

\noindent
{\it $q$-Lie algebras.} \
The discussion and examples in section 2 allows us to generalize and formulate
a general object as follows.

\proclaim{Definitions 3.1} For a $\Bbb Z$-graded vector space
$\Cal L=\bigoplus_{i\in\Bbb Z}\Cal L_i$ over a field $\Cal K$
equipped with a bilinear $q$-bracket product
$\{,\}_q$ $($where $q\in \Cal K$, $q\ne 0,1$,
$\dim\Cal L_i<\infty\,)$ satisfying
$\{\Cal L_i,\Cal L_j\}_q\subseteq \Cal L_{i+j}$, let
$\Cal S=\bigoplus_{i\in\Bbb Z}\Cal S_i\subseteq \Cal L$
be a  $\Bbb Z$-graded subspace of $\Cal L$.

\vskip0.2cm
\noindent
$(\text{\rm i})$ If the antisymmetry identity $(6)$ and the {\bf weighted $q$-Jacobi identity}
$(7)$ are fulfilled under $\{,\}_q$ for $x_i\in \Cal L_i,
\forall i\in \Bbb Z$, then $\Cal L^q:=(\Cal L,\{,\}_q)$ is called a
{\bf $q$-Lie algebra} and $\{,\}_q$ is called the {\bf $q$-Lie product}.
In particular, $(\Cal L_0,\{,\}_q|_{\Cal L_0})$ is a {\bf usual Lie algebra}
and $\{,\}_q|_{\Cal L_0}$ is a {\bf usual Lie product}.

\vskip0.2cm
\noindent
$(\text{\rm ii})$ If $\Cal S$ is closed under the $q$-Lie product $\{,\}_q$, then
$\Cal S^q$ is called a {\bf $q$-Lie subalgebra} of $\Cal L^q$. Let
$\Cal S\subset \Cal L$ be a subset of $\Cal L$, define a {\bf $q$-centralizer}
$\Cal C_q(\Cal S):=\{x\in \Cal L\mid \{x,\Cal S\}_q=0\}$ of $\Cal S$ in
$\Cal L$, and a {\bf $q$-normalizer} $\Cal N_q(\Cal S):=\{x\in\Cal L\mid
\{x,\Cal S\}_q\subseteq\Cal S\}$ of $\Cal S$ in $\Cal L$, then
both are $q$-Lie subalgebras, thanks to weighted $q$-Jacobi identities.
If moreover, a $q$-Lie subalgebra $\Cal S$ satisfies
$\{\Cal S,\Cal L\}_q\subseteq\Cal S$, then $\Cal S^q$ is called a
{\bf $q$-Lie ideal} of $\Cal L^q$. In particular, if
$\{\Cal S,\Cal L\}_q=0$, $\Cal S^q$ is called a {\bf $q$-central
subalgebra} of $\Cal L^q$. If $\{\Cal L,\Cal L\}_q=0$, then $\Cal L^q$ is
called {\bf $q$-abelian}. If $\Cal L^q$ only has trivial $q$-Lie ideals
$0$ and $\Cal L^q$, then $\Cal L^q$ is said {\bf $q$-simple}. For instance,
$\W^q$, $\W^q(1)$ and $\W^\vn(1,\bold 1)$. If $\Cal L^q$ is a finite direct
sum of $q$-simple ideals, then $\Cal L^q$ is {\bf $q$-semisimple}.

\vskip0.2cm
\noindent
$(\text{\rm iii})$ For any two $q$-Lie algebras $\Cal L^q$ and ${\Cal L'}^q$,
if there exists a graded vector space linear mapping
$\phi: \Cal L^q\lra{\Cal L'}^q$ of degree $0$ such that
$\phi(\Cal L_i)\subseteq\Cal L_i'$,
$\phi(\{x_i,x_j\}_q)=\{\phi(x_i),\phi(x_j)\}_q$, $\forall x_i\in\Cal L_i$
and $x_j\in\Cal L_j$,
then $\phi$ is called a {\bf $q$-Lie homomorphism}.
Then the {\bf kernel} of $\phi$ is a $q$-Lie ideal of $\Cal L^q$. The
notions of the {\bf $q$-Lie quotient algebra} and the {\bf image} of a
$q$-Lie homomorphism
also can be defined. The homomorphic fundamental theorem holds for
$q$-Lie algebras.

\vskip0.2cm
\noindent
$(\text{\rm iv})$ Given three $q$-Lie algebras $\Cal L^q$, ${\Cal L'}^q$ and
$\b {\Cal L}^q$, if ${\Cal L'}^q$ is a $q$-Lie subalgebra of $\b {\Cal L}^q$
and there exists an exact $($non-split$)$ sequence
$$
0\lra {\Cal L'}^q\lra \b{\Cal L}^q\lra \Cal L^q\lra 0,\eqno(16)
$$
then the $\b{\Cal L}^q$ is referred to as a $($nontrivial$)$ {\bf $q$-extension} of
$\Cal L^q$ through ${\Cal L'}^q$. In particular, if ${\Cal L'}^q$ is a
$q$-central subalgebra  of  $\b{\Cal L}^q$, then $\b{\Cal L}^q$ is said
a {\bf $q$-central extension} of $\Cal L^q$.

\vskip0.2cm
\noindent
$(\text{\rm v})$ Let $\Cal L^q$ be a $q$-Lie algebra, a vector space $\Cal M$ over
$\Cal K$ is called a {\bf left $\Cal L^q$-module}, if there exists a bilinear
map $\Cal L^q\times \Cal M\lra \Cal M$ such that
$$
\{x_i,x_j\}_q.m=q^{i+1}x_i.(x_j.m)-q^{j+1}x_j.(x_i.m),\eqno(17)
$$
for all $x_i\in\Cal L_i, x_j\in\Cal L_j, m\in\Cal M$. For instance,
$\Cal K[x,x^{-1}]$, $\frak A^q(1)$ and $\frak A^q(1,\bold 1)$ are
modules of $\W^q$, $\W^q(1)$ and $\W^\vn(1,\bold 1)$, respectively.

\vskip0.2cm
\noindent
$(\text{\rm vi})$ The associative algebra $\bold U(\Cal L^q):=\frak T(\Cal L^q)/\Cal J_q$ is
termed the {\bf (universal) enveloping algebra} of $\Cal L^q$,
where $\frak T(\Cal L^q)$ be the tensor algebra of $\Cal L^q$,
$\Cal J_q$ be the 2-sided ideal of $\frak T(\Cal L^q)$ generated by the
elements
$$
J(x_i,x_j):=q^{i+1}x_i\ot x_j-q^{j+1}x_j\ot x_i-\{x_i,x_j\}_q,
\quad\forall x_i\in\Cal L_i.\eqno(18)
$$
The composite map $\sigma$ of the canonical maps
$\Cal L^q\lra\frak T(\Cal L^q)\lra\bold U(\Cal L^q)$ is called the
{\it canonical map} of $\Cal L^q$ into $\bold U(\Cal L^q)$, then
in $\bold U(\Cal L^q)$, we have $q^{i+1}\sigma(x_i)\sigma(x_j)-q^{j+1}
\sigma(x_j)\sigma(x_i)=\sigma(\{x_i,x_j\}_q), \forall x_i, x_j$ $\in \Cal L^q$.

If $\Cal L^q$ is $q$-abelian, then $\Cal J_q=\Cal I_q$, $\bold U(\Cal L^q)$
is a quasipolynomial algebra $($cf. \cite{5}$)$, if moreover, $\dim \Cal L^q=n$,
then $\bold U(\Cal L^q)=\Cal K_q^n$ is a quantum $n$-space $($see section 2$)$.
\endproclaim

$\bold U(\Cal L^q)$ is {\bf universal} in the following sense.
\proclaim{Lemma 3.2}
Let $\sigma$ be the canonical map of $\Cal L^q$ into
$\bold U(\Cal L^q)$, let $\Cal A$ be an algebra with $1$, and let $\theta$
be a linear map of $\Cal L^q$ into $\Cal A$ such that
$$
q^{i+1}\theta(x_i)\theta(x_j)-q^{j+1}\theta(x_j)\theta(x_i)
=\theta(\{x_i,x_j\}_q), \quad \forall x_i\in\Cal L_i.
$$
Then there exists a unique homomorphism $\theta'$ of $\bold U(\Cal L^q)$ into
$\Cal A$ such that $\theta'(1)=1$ and $\theta' \cdot \sigma=\theta$.
\endproclaim
\demo{Proof} Since $\bold U(\Cal L^q)$ is generated by $1$ and
$\sigma(\Cal L^q)$, $\theta'$ is unique. On the other hand, let $\phi$
be the unique homomorphism of $\frak T(\Cal L^q)$ into $\Cal A$ which extends
$\theta$ with $\phi(1)=1$. $\forall x_i\in\Cal L_i$, we have
$$\split
&\phi(q^{i+1}x_i\ot x_j-q^{j+1}x_j\ot x_i-\{x_i,x_j\}_q)\\
&\quad=q^{i+1}\theta(x_i)\theta(x_j)
-q^{j+1}\theta(x_j)\theta(x_i)-\theta(\{x_i,x_j\}_q)\\
&\quad=0,
\endsplit
$$
hence $\phi(\Cal J_q)=0$ and, by passage to the quotient, $\phi$ defines a
homomorphism $\theta'$ of $\bold U(\Cal L^q)$ into $\Cal A$ with $\theta'(1)=1$
and $\theta'\cdot \sigma=\theta$.
\hfill\qed
\enddemo

\noindent
{\it Remark.} \
The bilinearity of $q$-Lie product $\{,\}_q|_{\Cal L^q}$ in
$\bold U(\Cal L^q)$ can be understood in the following sense:
for any $x=\sum x_i$, $x'=\sum x_j'\in \Cal L^q$, denote
$\tilde x:=\sum q^{i+1}x_i$, $\tilde x':=\sum q^{j+1}x_j'\in \Cal L^q$,
we define
$$\{x,x'\}_q:=\tilde xx'-\tilde x'x,\eqno(19)$$
then $\{x,x'\}_q=\{\sum x_i,\sum x_j'\}_q=\sum_{i,j}\{x_i,x_j'\}_q$.
So the ideal $\Cal J_q$ is generated by all $J(x,y)=
\tilde xy-\tilde yx-\{x,y\}_q$,
$\forall x, y\in\Cal L^q$.

\vskip0.3cm
\noindent
{\it $q$-PBW theorem.} \
Berger (\cite{2}) has studied a wide class of associative algebras $\U$
({\it affine} or {\it quadratic}) defined by generators and relations in which the
Poincar\'e-Birkhoff-Witt theorem is valid. This class includes numerious
recently appeared quantum algebras. Roughly speaking, our $\U(\Cal L^q)$
belongs to a class of {\it affine} $q$-algebras in the sense of Berger (see
section 2.1 in \cite{2}). But his proof applies only to the {\it quadratic} case
(see his Definition 2.5.2 and the proof of Lemma 2.8.2). According to his
remark 2.8.3, the proof of the {\it affine} case essentially follows
Jacobson's (cf. \cite{8}). We will give a proof for $\U(\Cal L^q)$ by following
Bergman's idea of Diamond Lemma --- reduction system (\cite{3}).

\vskip0.16cm
We need some notions and notations from \cite{3}.

\vskip0.16cm
Fix an ordered homogeneous basis $\{x_{i_1},x_{i_2},\cdots\}$
$(i_1< i_2<\cdots)$ of $\Cal L^q$ such that $x_{i_k}\in\Cal L_n$ for some $n$.
Let $\{\Cal L^q\}$ denote the set of monomials in $\frak T(\Cal L^q)$.
Define the {\it disordering index} of a monomial
$x_{j_1}\ot\cdots\ot x_{j_n}$ as the number of pairs $(j_l,j_k)$ such
that $j_l<j_k$ but $x_{j_l}>x_{j_k}$ (cf. \cite{8}). Thus we can partially
order monomials in $\{\Cal L^q\}$ by setting $\Cal C\prec\Cal D$ if
$\Cal C$ is of smaller length than $\Cal D$, or if $\Cal C$ is a permutation
of the terms of $\Cal D$ but of smaller disordering index. Clearly, this
defines a {\it semigroup partial ordering} on $\{\Cal L^q\}$ in the sense:
$\Cal B\prec\Cal B'$ implies
$\Cal A\ot\Cal B\ot\Cal C\prec\Cal A\ot\Cal B'\ot\Cal C$ $(\Cal {A, B, B', C}
\in \{\Cal L^q\})$.

\vskip0.16cm
Moreover, let $S=\{\,\tau=(W_{\tau},f_{\tau})\mid W_{\tau}\in\{\Cal L^q\},
f_\tau\in\frak T(\Cal L^q)\,\}$ be a {\it reduction system} of $\frak
T(\Cal L^q)$. If $\forall \tau\in S$, $f_\tau$ is
a linear combination of monomials $\prec W_\tau$, then the semigroup partial
ordering $\prec$ is called {\it compatible} with $S$.
For each $\tau\in S$, let $r_\tau$ denote the map on $\frak T(\Cal L^q)$ which
maps each monomial $\Cal A\ot W_\tau\ot\Cal B$ into $\Cal A\ot f_\tau\ot
\Cal B$ and fixes those monomials without containing subword $W_\tau$.
Denote $r_S:=\{\,r_\tau\mid \tau\in S\,\}$ and each $r_\tau$ is called a {\it
reduction} of $\frak T(\Cal L_q)$. For any $\Cal A\in\{\Cal L^q\}$, define
$\Cal J_{\Cal A}:=\langle \,\Cal C=\sum\Cal C_i\in\Cal J_q\mid \Cal C_i\prec
\Cal A, \text{\it for each } i\,\rangle$.
An {\it overlap ambiguity} of $S$ means a $5$-tuple $(\sigma,\tau,\Cal A,\Cal B,
\Cal C)$ with $\sigma,\tau\in S$ and $\Cal {A, B, C}\in \{\Cal L^q\}-{1}$ such
that $W_\sigma=\Cal A\ot\Cal B, W_\tau=\Cal B\ot\Cal C$. The ambiguity is said
{\it resovable relative to} $\prec$ if $f_\sigma\ot\Cal C-\Cal A\ot f_\tau\in
\Cal J_{\Cal A\ot\Cal B\ot\Cal C}$.

\vskip0.16cm
For any $x_i\in\Cal L_i, x_j\in\Cal L_j, x_k\in\Cal L_k$ such that $x_i<x_j<x_k$, the
{\it $q$-Jacobi sum} is the following element of ${\Cal L^q}^{\ot 2}$:
$$
\split
J(x_k,x_j,x_ i)&:=q^{j+k}\{x_k,x_j\}_q\ot x_i-q^{2i}x_i\ot \{x_k,x_j\}_q\\
&\quad +q^{i+j}\{x_j,x_i\}_q\ot x_k-q^{2k}x_k\ot \{x_j,x_i\}_q\\
&\quad +q^{i+k}\{x_i,x_k\}_q\ot x_j-q^{2j}x_j\ot \{x_i,x_k\}_q.
\endsplit\tag{20}
$$

The following lemma is essentially due to Berger (cf. Proposition 2.4.1
\cite{2}).
\proclaim{Lemma 3.3} The $q$-Jacobi sums $J(x_k,x_j,x_i)\in
\Cal J_{x_k\ot x_j\ot x_i}$.
\endproclaim
\demo{Proof}
Note that each monomial (in $J(x_k,x_j,x_i)$) $\prec x_k\ot x_j\ot x_i$.
We have to show $J(x_k,x_j,x_i)\in \Cal J_q$.

From (18), we have
$\{x_i,x_j\}_q=q^{i+1}x_i\ot x_j-q^{j+1}x_j\ot x_i-J(x_i,x_j)$,
then
$$
\split
J(x_k,x_j,x_i)&=-q^{j+k}J(x_k,x_j)\ot x_i+q^{2i}x_i\ot J(x_k,x_j)\\
&\quad -q^{i+j}J(x_j,x_i)\ot x_k+q^{2k}x_k\ot J(x_j,x_i)\\
&\quad -q^{i+k}J(x_i,x_k)\ot x_j+q^{2j}x_j\ot J(x_i,x_k)\in \Cal J_q.
\endsplit
$$

The proof is completed.
\hfill\qed
\enddemo

Denote by $\bar a$ the class in $\U(\Cal L^q)$ of any element $a$ of
$\frak T(\Cal L^q)$.
\proclaim{Theorem 3.4}
Let $\{x_{i_1},x_{i_2},\cdots\}$ $(i_1< i_2<\cdots)$ be an ordered homogeneous basis of $\Cal L^q$
 over $\Cal K$. Then the unit $1$ and the ordered monomials
$\bar x_{j_1}\cdots\bar x_{j_n}$ $(j_1\le\cdots\le j_n, n\ge 1)$
form a basis of $\U(\Cal L^q)$ over $\Cal K$.
\endproclaim
\demo{Proof}
Take $S=\{\,\tau=\tau_{xy}:=(W_{\tau_{xy}},f_{\tau_{xy}})\mid x<y,
x\in\Cal L_i, y\in\Cal L_j\,\}$ as a {\it reduction system} of $\frak T(\Cal L^q)$,
where $W_{\tau_{xy}}=q^{j+1}y\ot x\in\{\Cal L^q\}, f_{\tau_{xy}}=q^{i+1}x\ot y
-\{x,y\}_q\in\frak T(\Cal L^q)$ and the semigroup partial ordering $\prec$ is
compatible with $S$. Note that the ideal generated by the differences
$W_\tau-f_\tau$ $(\tau\in S)$ is precisely $\Cal J_q$ (since $\{,\}_q$ is
antisymmetry and bilinear!).
Since each reduction $r_\tau$
maps each monomial $\Cal A\ot W_\tau\ot\Cal B$ into $\Cal A\ot f_\tau\ot\Cal B$
but fixes those monomials without containing subword $W_\tau$, the images in
$\U(\Cal L^q)$ of the fixed monomials under $r_S$ are precisely the alleged
basis.

\vskip0.1cm
Since each monomial of length $n$ via certain finite number of reductions
will be stable under $r_S$, the above semigroup partial ordering $\prec$ on
$\{\Cal L^q\}$ satisfies the descending chain condition.

\vskip0.1cm
Clearly, the ambiguities of $S$ are precisely the $5$-tuples
$(\tau_{x_jx_k},\tau_{x_ix_j},$ $x_k,x_j,x_i)$ with $x_i<x_j<x_k$.
To see that the ambiguity is resovable relative to $\prec$, we study
$$
\split
&r_{\tau_{x_jx_k}}(q^{2k+j+1}x_k\ot x_j\ot x_i)
-r_{\tau_{x_ix_j}}(q^{2k+j+1}x_k\ot x_j\ot x_i)\\
&\quad=(q^{k+2j+1}x_j\ot x_k\ot x_i+q^{j+k}\{x_k,x_j\}_q\ot x_i)\\
&\qquad-(q^{2k+i+1}x_k\ot x_i\ot x_j+q^{2k}x_k\ot \{x_j,x_i\}_q).
\endsplit
$$
To further reduce the term $q^{k+2j+1}x_j\ot x_k\ot x_i$, we apply first
$r_{\tau_{x_ix_k}}$, and then $r_{\tau_{x_ix_j}}$. Similarly, to deal with
$q^{2k+i+1}x_k\ot x_i\ot x_j$ we apply $r_{\tau_{x_ix_k}}$ and then
$r_{\tau_{x_jx_k}}$. Thus we get
$$\split
&(q^{2i+j+1}x_i\ot x_j\ot x_k
+q^{i+j}\{x_j,x_i\}_q\ot x_k\\
&\quad+q^{2j}x_j\ot\{x_k,x_i\}_q
+q^{j+k}\{x_k,x_j\}_q\ot x_i)\\
&\  -(q^{2i+j+1}x_i\ot x_j\ot x_k
+q^{2i}x_i\ot \{x_k,x_j\}_q\\
&\quad+q^{i+k}\{x_k,x_i\}_q\ot x_j
+q^{2k}x_k\ot \{x_j,x_i\}_q)\\
&=J(x_k,x_j,x_i)\in\Cal J_{x_k\ot x_j\ot x_i}.
\endsplit
$$
So our ambiguities are resovable relative to $\prec$. Hence by Theorem 1.2
of Bergman in \cite{3}, $\U(\Cal L^q)$ has the basis indicated.
\hfill\qed
\enddemo

\proclaim{Corollary 3.5}
$(\text{\rm i})$ \ The canonical map $\sigma: \Cal L^q\lra \U(\Cal L^q)$ is injective.

$(\text{\rm ii})$ \ $\U(\Cal L^q)$ has no zero divisors.\hfill\qed
\endproclaim

\vskip0.3cm
\noindent
{\it $q$-Virasoro algebra $\text{\bv}^q$.} \
Here is another example of $q$-Lie algebras --- $q$-Virasoro algebra.
As an application of weighted $q$-Jacobi identity (7), we can use it to determine
the structure of the universal $q$-central extension of the $q$-Witt algebra
$\W^q$ defined in Lemma 2.1. A standard method of computing the central
extension in Lie theory (cf. \cite{1}, \cite{7} etc.) also applies to $q$-Lie
algebras here.

\vskip0.16cm
Recall that the Virasoro algebra $\text{\bv}=\langle e_i, c\mid i\in \Bbb Z\rangle$
is an universal $1$-dimensional central extension  of the Witt algebra
$\W=\langle e_i\mid i \in \Bbb Z\rangle$ with the Lie product
$$\gather
[e_i,c]=0,\\
[e_i,e_j]=(j-i)e_{i+j}+\delta_{i+j,0}\frac{(i-1)i(i+1)}{12}c.
\endgather
$$

Define $\text{\bv}^q=\langle e_i, L_0\mid i\in\Bbb Z\rangle$ with the
$q$-Lie product
$$\gather
\{e_i, L_0\}_q=0,\tag {21}\\
\{e_i,e_j\}_q=[(j+1)_q-(i+1)_q]\,e_{i+j}+
\delta_{i+j,0}\frac{(i-1)_q(i)_q(i+1)_q}{q^i(2)_{q^i}(2)_q(3)_q}\,L_0.\tag{22}
\endgather
$$

\proclaim{Proposition 3.6}
$\text{\bv}^q$ is a $q$-Lie algebra, which is an universal
$1$-dimensional $q$-central extension of the
$q$-Witt algebra $\W^q$, i.e. there is the exact sequence
$$0\lra \langle L_0\rangle\lra \text{\bv}^q\lra \W^q\lra 0.$$
\endproclaim
\demo{Proof} \ In this proof, we will derive the defining relation (22) of
$\text{\bv}^q$ using the weighted $q$-Jacobi identity.

Suppose that $0\lra \Cal C^q\lra \Cal L^q\lra \W^q\lra 0$ is a nonsplit
$q$-central extension of $\W^q$. Let $\tilde e_i$ denote the preimages of
$e_i$. Then we have
$$
\{\tilde e_i,\tilde e_j\}_q=o(i,j)\tilde e_{i+j}+c(i,j),\eqno(23)
$$
where $o(i,j)=[(j+1)_q-(i+1)_q]$, $c(i,j)=-c(j,i)$ and $c(i,j)\in\Cal C^q$.
Since $\{\tilde e_0,\tilde e_i\}_q=o(0,i)\tilde e_i+c(0,i)$. Replace
$\tilde e_i$ by $\tilde e_i+\frac{c(0,i)}{o(0,i)}$,
then we have
$$
\{\tilde e_0,\tilde e_i\}_q=o(0,i)\tilde e_i.
$$
For $\{\tilde e_1,\tilde e_{-1}\}_q=o(1,-1)\tilde e_0+c(1,-1)$, we also
can take $\tilde e_0$ such that $c(1,-1)=0$, i.e.
$$
\{\tilde e_1,\tilde e_{-1}\}_q=o(1,-1)\tilde e_0.
$$

First, we will show that $c(i,j)=0$,
for $i+j\ne 0$ ($\forall i,j\in\Bbb Z$).
From the following weighted $q$-Jacobi identity:
$$
\split
&\quad=(2)_{q^i}o(0,j)\{\tilde e_i,\tilde e_j\}_q\\
&(2)_{q^i}\{\tilde e_i,\{\tilde e_0,\tilde e_j\}_q\}_q\\
&\quad=(2)_{q^0}\{\tilde e_0,\{\tilde e_i,\tilde e_j\}_q\}_q+(2)_{q^j}\{\{\tilde e_i,\tilde e_0\}_q,\tilde e_j\}_q\\
&\quad=2o(i,j)o(0,i+j)\tilde e_{i+j}-(2)_{q^j}o(0,i)\{\tilde e_i,\tilde e_j\}_q,
\endsplit
$$
we get $[(2)_{q^i}o(0,j)+(2)_{q^j}o(0,i)]c(i,j)=2o(0,i+j)c(i,j)=0$, i.e.
$c(i,j)=0$, if $i+j\ne 0$.

\vskip0.16cm
Next, we will determine $c(i,-i)$.
For $r=i+j\ge 3$, using the weighted $q$-Jacobi identity:
$$\split
&\quad=(2)_{q^r}o(-i,-j)\bigl[o(r,-r)\tilde e_0+c(r,-r)\bigr]\\
&(2)_{q^r}\{\tilde e_r,\{\tilde e_{-i},\tilde e_{-j}\}_q\}_q\\
&\quad=(2)_{q^{-j}}\{\{\tilde e_r,\tilde e_{-i}\}_q,\tilde e_{-j}\}_q
      +(2)_{q^{-i}}\{\tilde e_{-i},\{\tilde e_r,\tilde e_{-j}\}_q\}_q\\
&\quad=(2)_{q^{-j}}o(r,-i)\bigl[o(j,-j)\tilde e_0+c(j,-j)\bigr]\\
&\quad\,-(2)_{q^{-i}}o(r,-j)\bigl[o(i,-i)\tilde e_0+c(i,-i)\bigr],
\endsplit
$$
we obtain
$$\split
&(2)_{q^r}o(-i,-j)c(r,-r)\\
&\quad =(2)_{q^{-j}}o(r,-i)c(j,-j)
-(2)_{q^{-i}}o(r,-j)c(i,-i).
\endsplit\tag{24}
$$
Note that $(2)_{q^{-j}}=q^{-j}(2)_{q^j}$,
and set $\Delta(r)=q^r(2)_{q^r}c(r,-r)$,  then
(24) becomes
$$
\Delta(r)(q^{-i}-q^{-j})
=\Delta(j)(q^{2i}-q^{-j})-\Delta(i)(q^{2j}-q^{-i}).$$
Since $c(1,-1)=0$, i.e. $\Delta(1)=0$, take $i=1, j=r-1$, then
we get
$$
\Delta(r)=\frac{q^2-q^{1-r}}{q^{-1}-q^{1-r}}\Delta(r-1)
=\frac{(r+1)_q}{(r-2)_q}\Delta(r-1).\eqno(25)
$$
It is easy from (25) to get
$$
\Delta(i)=\frac{(i+1)_q(i)_q(i-1)_q}{(2)_q(3)_q}\Delta(2).
$$
Set $L_0:=\Delta(2)$, we obtain the required conclusion:
$$
c(i,-i)=\frac{(i-1)_q(i)_q(i+1)_q}{q^i(2)_{q^i}(2)_q(3)_q}\,L_0.\eqno(26)
$$
Hence, $\Cal C^q=\langle L_0 \rangle$ and $\Cal L^q=\text{\bv}^q$.

\vskip0.16cm
Finally, we shall show that the product $\{,\}_q$ defined on $\text{\bv}^q$
(cf. (21), (22)) is  a $q$-Lie product.

\vskip0.1cm
At first, it is not hard to
see that $c(-i,i)=-c(i,-i)$ using the fact $(-i)_q=-q^{-i}(i)_q$, namely,
the antisymmetry holds for the $\{,\}_q$.

\vskip0.1cm
It remains to verify that the $\{,\}_q$
satisfies the weighted Jacobi identity. In view of $\text{\bv}^q$ being a
$1$-dimensional $q$-central extension of the $q$-Lie algebra $\W^q$, we
only need to check it for the case when $r=-(i+j)\not=0$, namely, to show
the following identity:
$$
(2)_{q^i}\{\t e_i,\{\t e_j,\t e_r\}_q\}_q+
(2)_{q^j}\{\t e_j,\{\t e_r,\t e_i\}_q\}_q+
(2)_{q^r}\{\t e_r,\{\t e_i,\t e_j\}_q\}_q=0.\eqno(27)
$$
Applying (22) to the left side of (27), it suffices to verify that
the sum  of the coefficents of
$L_0$ in (27),
$\Cal S=\Cal S(i,j,r)+\Cal S(j,r,i)+\Cal S(r,i,j)$ is zero (since the fact
that the sum of the coefficents of $\t e_0$ in (27) is zero is due to the
weighted Jacobi identity of $\W^q$).

For $k\in\Bbb Z$, denote $[k]:=q^k+q^{-k}$, then
$[-k]=[k]$. It is easy to see that the coefficient of $L_0$ in
$(2)_{q^i}\{\t e_i,\{\t e_j,\t e_r\}_q\}_q$ is
$$\split
\Cal S(i,j,r)&=\frac{1}{(2)_q(3)_q}q^{-i}(i-1)_q(i)_q(i+1)_q\bigl((r+1)_q-(j+1)_q\bigr)\\
&=\frac{q}{(2)_q(3)_q}\bigl([i-j]-[i-r]-q^{-1}(3)_q[j]+q^{-1}(3)_q[r]
\bigr).
\endsplit
$$
Thus $\Cal S=0$. This shows the $\{,\}_q$ is indeed a $q$-Lie product and
$\text{\bv}^q$ is a $q$-Lie algebra.

The proof is completed.
\hfill\qed
\enddemo

\vskip0.25cm
\noindent
{\bf 4 \ Algebra $\bold U^\vn$}

\vskip0.3cm
\noindent
From now on, our interest mainly concentrates on
the situation of primitive $l$-th root of unity.
Some of particular observations only can be carried out in this case but
the results obtained in this way can be transferred to the generic case.

\vskip0.2cm
\noindent
{\it Algebra\ $\bold U^\vn$.} \
For the $\vn$-Witt algebra $\W^\vn(1,\bold 1)$ over $\Cal K$, consider
its universal enveloping algebra $\bold U^\vn:=\bold U(\W^\vn(1,\bold 1))$,
which contains a big central subalgebra analogous to the modular case.

Set $H_{ij}^{(0)}=1$, and for $k>0$,
$$
H_{ij}^{(k)}
=\cases
\vn^{{-}k{-}j\binom{k{+}1}{2}}\prod_{s=1}^k
\left[\binom{i{+}sj{+}1}{j}_\vn{-}\binom{i{+}sj{+}1}{j{+}1}_\vn\right],
 & {-}1\le i{+}j\le l{-}2,\\
0, & \text{\it otherwise}.
\endcases
$$

\proclaim
{Lemma 4.1} \ {\it For } $-1\le i, j\le l-2$, {\it we have}

\vskip0.08cm
\;\;$(\text{\rm i})$ \;\quad $e_{(i)}\,{e_{(0)}}^n=\vn^{-in}(e_{(0)}-(i)_\vn)^ne_{(i)}$;

\vskip0.12cm
$(\text{\rm ii})$ \quad $e_{(i)}\,{e_{(j)}}^n
      =\sum_{k=0}^n\binom{n}k_{\vn^{-j}}
      \vn^{(n-k)(j-i)}{e_{(j)}}^{n-k}[\cdots [e_{(i)},
      \undersetbrace k\roman{\;times}\to{e_{(j)}]_\vn,\cdots,e_{(j)}}]_\vn$

\hskip2.34cm $=\vn^{n(j-i)}\sum_{k=0}^n\binom{n}k_{\vn^{-j}}
      H_{ij}^{(k)}{e_{(j)}}^{n-k}e_{(i+kj)}$,\quad $(j\not=0)$.
\endproclaim
\demo{Proof} (i) follows from $e_{(i)}e_{(0)}=\vn^{-i}(e_{(0)}-(i)_\vn)e_{(i)}$,
using (13).

(ii) Denote $[e_{(i)},e_{(j)}]_\vn:=e_{(i)}e_{(j)}-\vn^{j-i}e_{(j)}e_{(i)}$.
Assume that the first equality of (ii) holds for $n$. Since by (13),
$[e_{(i)},e_{(j)}]_\vn=*\;e_{(i+j)}$, then
$$\split
[\cdots[e_{(i)},\undersetbrace k\to{e_{(j)}]_\vn,\cdots,e_{(j)}}]_\vn e_{(j)}
&=\vn^{j-(i+kj)}e_{(j)}[\cdots [e_{(i)},
\undersetbrace k\to{e_{(j)}]_\vn,\cdots,e_{(j)}}]_\vn\\
&\qquad+[\cdots [e_{(i)},\undersetbrace {k+1}\to{e_{(j)}]_\vn,\cdots,e_{(j)}}]_\vn,\\
\endsplit
$$
and
$$\split
e_{(i)}{e_{(j)}}^{n+1}
&=
\sum_{k=0}^n\binom{n}{k}_{\vn^{{-}j}}
\vn^{(n{-}k)(j{-}i){+}j{-}(i{+}kj)}{e_{(j)}}^{n{+}1{-}k}
[{\cdots}[e_{(i)},\undersetbrace k\to{e_{(j)}]_\vn,{\cdots},e_{(j)}}]_\vn\\
&\qquad+
\sum_{k=0}^n\binom{n}{k}_{\vn^{-j}}
\vn^{(n-k)(j-i)}{e_{(j)}}^{n-k}
[\cdots[e_{(i)},\undersetbrace {k+1}\to{e_{(j)}]_\vn,\cdots,e_{(j)}}]_\vn\\
&\ =\sum_{k=0}^n\vn^{-kj}\binom{n}{k}_{\vn^{-j}}
\vn^{(n+1-k)(j-i)}{e_{(j)}}^{n+1-k}
[\cdots[e_{(i)},\undersetbrace k\to{e_{(j)}]_\vn,\cdots,e_{(j)}}]_\vn\\
&\qquad+
\sum_{k=1}^{n+1}\binom{n}{k{-}1}_{\vn^{-j}}
\vn^{(n+1-k)(j-i)}{e_{(j)}}^{n+1-k}
[\cdots[e_{(i)},\undersetbrace k\to{e_{(j)}]_\vn,\cdots,e_{(j)}}]_\vn.\\
\endsplit
$$
But
$$
\binom{n}{k{-}1}_{\vn^{{-}j}}
+\vn^{-kj}\binom{n}{k}_{\vn^{{-}j}}
=\binom{n{+}1}{k}_{\vn^{{-}j}}.
$$
Hence, the first equality of (ii) holds for $n+1$.

The second equality of (ii) is obvious due to (13).
\hfill\qed
\enddemo

\proclaim
{Corollary 4.2} \ {\it For } $-1\le i,j\le l-2$, {\it we have}

\vskip0.08cm
$(\text{\rm i})$ \;\quad $\bigl[e_{(i)},({e_{(0)}-\frac{1}{1-\vn})}^{l}\bigr]=0$;

\vskip0.08cm
$(\text{\rm ii})$ \quad $\bigl[e_{(i)},{e_{(j)}}^{l}\bigr]=0$\quad{\it for }
\ $j\not =0$.
\endproclaim
\demo{Proof} (i)
Set $\tilde e_{(0)}:=e_{(0)}-\frac{1}{1-\vn}$,
since $e_{(0)}e_{(i)}=\vn^ie_{(i)}e_{(0)}+(i)_\vn e_{(i)}$, by (13),
then
$\tilde e_{(0)}.e_{(i)}=\vn^ie_{(i)}.\tilde e_{(0)}$, and
${e_{(i)}.{\tilde e}_{(0)}}^{l}{\;={\tilde e}_{(0)}}^{l}.e_{(i)}$.

(ii) is clear by virtue of Lemma 4.1 (ii) and the facts that $e_{(i+lj)}=0$
for $-1\le i\le l-2$ and $\binom{l}k_{\vn^{-j}}=0$ for $0<k<l$.
\hfill\qed
\enddemo

Set $z_i={e_i}^{l}$ $(i\not =0)$ and
$z_0=(e_{(0)}-\frac{1}{1-\vn})^{l}$.
Then
$\Cal Z_0:=\Cal K[z_{-1},z_0,\cdots,z_{l-2}]$
is a central subalgebra of $\bold U^\vn$.

\vskip0.15cm
Following \cite{5}, we introduce a $\Bbb Z_+^l$-gradation in
$\bold U^\vn$ by letting $\text{\it deg}\,(e_{(i)})$
$=(\delta_{i,-1},{\cdots},\delta_{i,l-2})$ $(-1\le i\le l-2)$. For
an integer vector $\bold n:=(n_{-1},{\cdots},n_{l-2})$, $n_i\in \Bbb Z_+$, we
set $\text{\it deg}\,\bold n:=\sum n_i$ and $e^{\bold n}=e_{-1}^{n_{-1}}\cdots
e_{l-2}^{n_{l-2}}$ and call such an element a {\it monomial}. Furthermore
we define on the set of integer vectors the degree-lexicographic ordering,
i.e. set $\bold n<\bold m$ if either $\text{\it deg}\,\bold n<\text{\it deg}\,
\bold m$ or $\text{\it deg}\,\bold n=\text{\it deg}\,\bold m$ but $\bold n$ is
less than $\bold m$ in the usual lexicographic order, in this way
$\Bbb Z_+^l$ becomes an ordered monoid.

Denote
$$
\bold U^\vn_{\bold n}:=\sum_{{\bold m}\le{\bold n}} \Cal K e^{\bold m},
$$
then the subspaces $\bold U^\vn_{\bold n}$ give a structure of filtered algebra
relative to the ordered monoid $\Bbb Z_+^{l}$. Obviously, with respect to the
ordered monoid $\Bbb Z_+^{l}$, the associated graded algebra
$\text{\it gr}\,\bold U^\vn$ to the filtered algebra $\bold U^\vn$ has the relations:
$\bar e_{(i)}\bar e_{(j)}=\lambda_{ij}\bar e_{(j)}\bar e_{(i)}\quad
(\lambda_{ij}=\vn^{j-i})$,
which is isomorphic to a quasipolynomial algebra $\Cal P$ (cf. \cite{5}),
i.e. $\text{\it gr}\,\bold U^\vn\cong \Cal P$.
It shows again that $\bold U^\vn$ has no zero divisors
since $\Cal P$ has no zero divisors (cf. \cite{5}).
\proclaim{Proposition 4.3} {\it $\bold U^\vn$ is a maximal order.}
\endproclaim
\demo{Proof} \
Since $\bold U^\vn$ has no zero divisors, one may consider $\bold U^\vn$ as a free
left $\Cal Z_0$-module, and let
$$
\bold U^\vn_{\bold n}:=\sum_{\bold m\le\bold n}\Cal Z_0e^{\bold m},
$$
then, similarly, $\bold U^\vn$ has a structure of filtered algebra with coefficients
in $\Cal Z_0$ relative to the ordered monoid $\Bbb Z_{(l)}^{l}$
($\Bbb Z_{(l)}=\Bbb Z/(l)$) and the associated
graded algebra $\text{\it gr}\,\bold U^\vn$ is a twisted polynomial ring over
$\Cal Z_0$. Since $\Cal Z_0$ is integrally closed and $\bold U^\vn$ satisfies
evidently hypotheses 1--6 in the section 6.5 \cite{5}, by Theorem 6.5 \cite{5},
we get the conclusion.
\hfill\qed
\enddemo

Let $\Cal Z$ be the center of $\bold U^\vn$,
since $\bold U^\vn$ has no zero divisors, the algebra $Q(\bold U^\vn):
=Q(\Cal Z)\ot_{\Cal Z} \bold U^\vn$ is a division algebra of dimension
$d^2$ over the field of fractions
$Q(\Cal Z)$ of $\Cal Z$ and $d$ is called the degree of $\bold U^\vn$.
Moreover, since $\bold U^\vn$ is a maximal order in
$Q(\bold U^\vn)$, then $\Cal Z$ is integrally closed.

\vskip0.2cm
Let $\Cal A$ be an algebra over $\Bbb C[q,q^{-1}]$ on generators $x_1,
\cdots,x_n$ satisfying the relations:
$x_ix_j=q^{h_{ij}}x_jx_i+\Cal P_{ij}$, for $i>j$,
$\Cal P_{ij}\in\Bbb C[q,q^{-1}][x_1,\cdots,x_{i-1}]$ and
$\bigl(h_{ij}\bigr)$ is a skew-symmetric matrix over $\Bbb Z$.
Let $l>1$ be an integer
relatively prime to all elementary divisors of the matrix $\bigl(h_{ij}\bigr)$
and let $\vn$ be a primitive $l$-th root of $1$. Let
$\Cal A_\vn=\Cal A/(q-\vn)$ and assume that ${x_i}^l$ are central. Let
$\Cal Z_0=\Bbb C[{x_1}^l,\cdots,{x_n}^l]$, which has a canonical
Poisson structure. C. de Concini and C. Procesi (\cite{5})
suggested the following conjecture.
\proclaim
{Conjecture 4.4} {\it Let $\pi$ be an irreducible representation of the algebra
$\Cal A_\vn$ and let $\Cal O_\pi\subset \text{Spec}\,\Cal Z_0$ be the
symplectic leaf containing the restriction of the central character of
$\pi$ to $\Cal Z_0$. Then the dimension of this representation is equal to
$l^{\frac{1}2\dim\Cal O_{\pi}}$}.
\endproclaim

Set
$\W^\vn(1,\bold 1)_{(i)}=\text{span}_\Cal K\{e_{(i)}\}$,
$\W^\vn(1,\bold 1)_i=\bigoplus_{j\ge i}\W^\vn(1,\bold 1)_{(j)}$.
So the $\vn$-Witt algebra
$\W^\vn(1,\bold 1)$ is graded with $\text{\it deg}\,(e_{(i)})=i$, i.e.
$\W^\vn(1,\bold 1)=\W^\vn(1,\bold 1)_{-1}=
\bigoplus_{j=-1}^{l-2}\W^\vn(1,\bold 1)_{(j)}$.
Let $\bold U_i^\vn$ $(i\ge 0)$ denote the subalgebras generated by
$\W^\vn(1,\bold 1)_i$, and
let $\Cal Z_0^{(i)}=\Cal Z_0\cap \bold U_i^\vn=\Cal K[z_i,\cdots,z_{l-2}]$,
then the algebras $\bold U_i^\vn$ are of the same type as $\Cal A_\vn$.
In particular, for $\bold U^\vn_1$, we have a similar conjecture concerning
its irreducible representations. However, the treatment of this problem is
much more complicated than that of the modular Witt algebra.

\vskip0.5cm
\noindent
{\bf 5 \ Realization}
\vskip0.3cm

\noindent
{\it Realization.} \
In this section, an operator realization of a class of
representations of $\W^\vn(1,\bold 1)$ over a tensor vector space $\frak
A^\vn(1,\bold 1)\ot_{\Cal K} V$ is constructed, which to some extent can be
compared with Shen's mixed product realization (cf. \cite{11}) in the modular
case.

A $\W^\vn(1{,}\bold 1)$-module $\Cal M=\bigoplus_{i=0}^t\Cal M_{(i)}$ is
graded, if $e_{(i)}.\Cal M_{(j)}$ $\subseteq\Cal M_{(i+j)}$, where
$\bold W^\vn(1,\bold 1)_{(0)}$-modules $\Cal M_{(0)}$ and $\Cal M_{(t)}$ are
called the base space and the top space of $\Cal M$, respectively.  Then we have
\proclaim
{Theorem 5.1} Let $\frak A^\vn(1,\bold 1)$ be an $\vn$-divided power algebra and
$\tau_\vn$ the automorphism of $\frak A^\vn(1,\bold 1)$ as in section 2.
Assume that $V$ is a $\W^\vn(1,\bold 1)_{(0)}$-module and $\rho_0$ is the
associated representation. Then
the following linear map
$$
e_{(i)}\longmapsto
x^{(i+1)}\partial_\vn\ot\text{\it id}_V+
\partial_\vn(x^{(i+1)})\tau_\vn\ot
\rho_0(e_{(0)})\eqno {(28)}
$$
gives a realization of graded representations of $\W^\vn(1,\bold 1)$
on the vector space $\frak A^\vn(1,\bold 1)\ot_{\Cal K} V$
and
$x^{(0)}\ot V\cong V$
as $\W^\vn(1,\bold 1)_{(0)}$-modules.
\endproclaim
\demo{Proof}
For $x^{(a)}\ot v\in\frak A^\vn(1,\bold 1)\ot_{\Cal K} V$ with
$\text{\it deg}\,(x^{(a)}\ot v)=a$, (28) gives
$$\split
\ e_{(i)}.x^{(a)}\ot v &=x^{(i+1)}\partial_\vn(x^{(a)})\ot v+
\partial_\vn(x^{(i+1)})\tau_\vn(x^{(a)})\ot e_{(0)}.v\\
&=\binom{a{+}i}{i{+}1}_\vn x^{(a{+}i)}\ot v+\vn^a\binom{a{+}i}{i}_\vn x^{(a{+}i)}
\ot e_{(0)}.v.
\endsplit\tag{29}
$$
Then
$$\split
e_{(i)}.&e_{(j)}. x^{(a)}\ot v =\binom{a{+}i{+}j}{i,j}_\vn \Bigl[
\frac{(a)_\vn(a{+}j)_\vn}{(i{+}1)_\vn(j{+}1)_\vn}x^{(a{+}i{+}j)}\ot v\Bigr.\\
&+\frac{(a{+}j)_\vn\vn^a}{(i{+}1)_\vn}x^{(a{+}i{+}j)}\ot e_{(0)}.v
+\frac{(a)_\vn\vn^{a{+}j}}{(j{+}1)_\vn}x^{(a{+}i{+}j)}\ot e_{(0)}.v\\
&+\left.\vn^{2a{+}j}x^{(a{+}i{+}j)}\ot e_{(0)}^2.v\right],
\endsplit
$$
where $\binom{a+i+j}{i,j}_\vn
=\frac{(a+i+j)_\vn!}{(a)_\vn!(i)_\vn!(j)_\vn!}$.

Thus
$$\split
(&\vn^{i+1}e_{(i)}.e_{(j)}-\vn^{j+1}e_{(j)}.e_{(i)}).x^{(a)}\ot v\\
&\ =\binom{a{+}i{+}j}{i,j}_\vn\frac{(j{+}1)_\vn-(i{+}1)_\vn}{(i{+}1)_\vn(j{+}1)_\vn}
\left[(a)_\vn x^{(a{+}i{+}j)}\ot v
+(i{+}j{+}1)_\vn\vn^a x^{(a{+}i{+}j)}\ot e_{(0)}.v\right]\\
&\ =\left[\binom{i{+}j{+}1}{i{+}1}_\vn-\binom{i{+}j{+}1}{j{+}1}_\vn\right]
\left[\binom{a{+}i{+}j}{i{+}j{+}1}_\vn x^{(a{+}i{+}j)}\ot v
+\vn^a \binom{a{+}i{+}j}{i{+}j}_\vn x^{(a{+}i{+}j)}\ot e_{(0)}.v\right]\\
&\ =\left[\binom{i{+}j{+}1}{i{+}1}_\vn-\binom{i{+}j{+}1}{j{+}1}_\vn\right]e_{(i{+}j)}.x^{(a)}\ot v.
\endsplit
$$
This shows that (28) does yield a representation of $\W^\vn(1,\bold 1)$, owing
to (29) and (13).\hfill\qed
\enddemo

Obviously, the proof of Theorem 5.1 also applies to the
generic case.
\proclaim{Corollary 5.2} Let $\frak A^q(1)$ be a $q$-divided power algebra and
$\tau_q$ the automorphism of $\frak A^q(1)$ as in section 2.
Assume that $V$ is a $\bold W^q(1)_{(0)}$-module
and $\rho_0$ is the associated representation. Then the following linear map
$$
e_{(i)}\longmapsto
x^{(i+1)}\partial_q\ot\text{\it id}_V+
\partial_q(x^{(i+1)})\tau_q\ot
\rho_0(e_{(0)})\eqno {(30)}
$$
gives a realization of graded representations of $\W^q(1)$
on the vector space $\frak A^q(1)\ot_{\Cal K} V$ and $x^{(0)}\ot V\cong V$
as $\bold W^q(1)_{(0)}$-modules.
\hfill\qed
\endproclaim

\vskip0.25cm
\noindent
{\it Irreducible modules.} \
Based on the realization above, the irreducible graded
$\W^\vn(1,\bold 1)$-modules can be explicitly described.
First of all, the following fact is clear.
\proclaim{Lemma 5.3}
{\it For every irreducible $\W^\vn(1,\bold 1)_{(0)}$-module $V$,
there exists, up to isomorphism, one and only one irreducible graded
$\W^\vn(1,\bold 1)$-module $\Cal M$ with $\Cal M_{(0)}=V$ as its base space.}
\hfill\qed
\endproclaim

Clearly, every irreducible $\W^\vn(1,\bold 1)_{(0)}$-module $V$
is $1$-dimensional,
i.e. $e_{(0)}.v=tv$, $t\in \Cal K$.
Denote by $\Cal M(t)$ the associated
irreducible graded $\W^\vn(1,\bold 1)$-module.
Note that the base (resp. top) space of
$\frak A^\vn(1,\bold 1)\ot_{\Cal K} V(t)$ is $V(t)$ (resp. $V((-1)_\vn(1-t)$),
according to the formula (29).

\proclaim{Theorem 5.4} \item{}
$(\text{\rm i})$ \ \;If \ $t\not=0$, $1$, then the graded
\, $\W^\vn(1,\bold 1)$-modules \,$\frak A^\vn(1,\bold 1)$ $\ot_{\Cal K} V(t)$ with
the base spaces $V(t)$

\item{}
\qquad are irreducible.

\item{}
$(\text{\rm ii})$ \ If \ $t=0$, then the graded
\; $\W^\vn(1,\bold 1)$-module \;$\frak A^\vn(1,\bold 1)\ot_{\Cal K} V(0)$ with the
base space \,$V(0)$ is

\item{}
\qquad reducible and there is an exact sequence
$$
0\lra \Cal M(0)\lra \frak A^\vn(1,\bold 1)\ot_{\Cal K} V(0)\lra \Cal M(1)\lra 0,
$$
\qquad where $\Cal M(0)=\langle x^{(0)}\ot v\rangle$, $\Cal M(1)=\langle
x^{(1)}\ot v, x^{(2)}\ot v,\cdots,x^{(l-1)}\ot v\rangle$.

\item{}
$(\text{\rm iii})$ \ If \ $t=1$, then the graded \; $\W^\vn(1,\bold 1)$-module
\;$\frak A^\vn(1,\bold 1)\ot_{\Cal K} V(1)$ with the base space $V(1)$ is

\item{}
\qquad reducible and there is an exact sequence
$$
0\lra \Cal M(1)\lra \frak A^\vn(1,\bold 1)\ot_{\Cal K} V(1)\lra \Cal M(0)\lra 0,
$$
\qquad where $\Cal M(0)=\langle x^{(l-1)}\ot v\rangle$, $\Cal M(1)=\langle
x^{(0)}\ot v, x^{(1)}\ot v,\cdots,x^{(l-2)}\ot v\rangle$.

\item{}
\qquad Note that the cases $(\text{\rm ii})$ and $(\text{\rm iii})$ are dual to each other and the
$\W^\vn(1,\bold 1)$-module $\Cal M(1)$ is

\item{}
\qquad self-dual.
\endproclaim
\demo{Proof} (i) Since $e_{(-1)}.\,x^{(0)}\ot v=0$,
$e_{(-1)}.\,x^{(m)}\ot v=x^{(m-1)}\ot v$, for $m>0$,
an arbitrary nontrivial $\W^\vn(1,\bold 1)$-submodule of $\frak A^\vn(1,\bold 1)\ot_{\Cal K} V$
must contain some nonzero element
of the form $x^{(0)}\ot v$, for $v\not =0\in V$.
By (29), we obtain
$e_{(l-2)}.\,x^{(0)}\ot v=t\,x^{(l-2)}\ot v$ and
$e_{(1)}.\,x^{(l-2)}\ot v=(1-t)\,\vn^{-3}x^{(l-1)}\ot v$.
So statement (i) holds.

(ii) Since $t=0$, $e_{(i)}.\,x^{(0)}\ot v=0$, $\forall \ i\ge 0$, by (29).
Then $\langle x^{(0)}\ot v\rangle$ is a trivial one-dimensional
$\W^\vn(1,\bold 1)$-submodule of $\frak A^\vn(1,\bold 1)\ot_{\Cal K} V$. And since
$e_{(l-2)}.\,x^{(1)}\ot v=x^{(l-1)}\ot v$ by (29), then
$\{x^{(1)}\ot v,\cdots,x^{(l-1)}\ot v\}$ is an
irreducible $\W^\vn(1,\bold 1)$ quotient module.
Since by (29),
$e_{(0)}.\,x^{(1)}\ot v=x^{(1)}\ot v$ and
$e_{(0)}.\,x^{(l-1)}\ot v=(-1)_\vn x^{(l-1)}\ot v$.
So statement (ii) is true.

(iii) Since $t=1$, by (29), we have
$e_{(l{-}2)}.\,x^{(0)}\ot v
 =x^{(l{-}2)}\ot v$,
$e_{(i)}.\,x^{(l{-}2)}\ot v
 =0$ for $i\ge 2$ and
$e_{(1)}.\,x^{(l{-}2)}\ot v =\binom{l{-}1}{2}_\vn x^{(l{-}1)}\ot v{+}
\vn^{{-}2}\binom{l{-}1}{1}_\vn x^{(l{-}1)}\ot v=0$.
On the other hand, from (29), we get
$e_{(0)}.\,x^{(0)}\ot v
  =x^{(0)}\ot v$,
$
e_{(0)}.\,x^{(l{-}2)}\ot v
 =\binom{l{-}2}{1}_\vn x^{(l{-}2)}\ot v{+}\vn^{-2}x^{(l{-}2)}\ot v
=(-1)_\vn x^{(l{-}2)}\ot v$
and $e_{(0)}.\,x^{(l{-}1)}\ot v
  =\binom{l{-}1}{1}_\vn x^{(l{-}1)}\ot v{+}
     \vn^{l{-}1}x^{(l{-}1)}\ot v=0$.
Then statement (iii) holds.
\hfill\qed
\enddemo

\noindent
{\it Remark.} \ The conclusions (ii) and (iii) in the above Theorem
show that there exists a certain duality between the modules
$\frak A^\vn(1,\bold 1)\ot_{\Cal K} V(0)$ and
$\frak A^\vn(1,\bold 1)\ot_{\Cal K} V(1)$.

\vskip0.5cm

\noindent
{\bf 6 \ $q$-Holomorph Structure and Representations}

\vskip0.3cm

\noindent
{\it An approach towards more representations.} \,
In (14) and (15), we constructed the $\vn$-holomorph structure $\Cal H^\vn(1)$ for
the $\vn$-Witt algebra $\W^\vn(1,\bold 1)$. Here we will reveal the importance
of such a structure in the representation theory.

For the $\vn$-Witt algebra $\W^\vn(1,\bold 1)$, $\vn$-abelian algebra
$_1\Cal K_\vn^l$ and their $q$-holomorph $\Cal H^\vn(1)$, given two
representation vector spaces $\frak V$ and $\Cal V$ with three associated
representations:
$$
\split
\phi: \ &\W^\vn(1,\bold 1)\lra \text{End}_{\Cal K}(\frak V),\\
\psi: \ & _1\Cal K_\vn^l\lra \text{End}_{\Cal K}(\frak V),\\
\rho: \ &\Cal H^\vn(1)\lra \text{End}_{\Cal K}(\Cal V).
\endsplit
$$
and the pair ($\phi, \psi$) satisfies the following compatibility condition:
$$
\phi(e_{(i)})\psi(L_j)-\vn^{j-i}\psi(L_j)\phi(e_{(i)})
=\binom{i+j}{i+1}_\vn\psi(L_{i+j}).\eqno(31)
$$
Then we have
\proclaim{Theorem 6.1}  Suppose that $\phi$, $\psi$ and $\rho$ are as
above, and set $L_{-1}=0$.

\item{}
$(\text{\rm i})$ \ \;If $\phi$ and $\psi$ satisfy $(31)$, then the pairs
$(\phi_{a\psi},\psi) (\forall a\in \Cal K)$ afford a family of representations

\item{}
\qquad
of $\Cal H^\vn(1)$, where $\phi_{a\psi}(e_{(i)})=\phi(e_{(i)})+a\psi(L_i)$
$(i\ge 0)$.

\item{}
$(\text{\rm ii})$ \ For any fixed element $\omega\in \;\Cal H^\vn(1)$ with $\rho(\omega)\in
\text{End}_{\Cal K}(\Cal V)$,
if $(31)$ holds for the pair $(\phi,\;\psi)$,

\item{}
\qquad then the triple $(\phi,\psi;\rho)$
provides a new representation $\tilde \phi(e_{(i)}):
=\phi(e_{(i)})\ot \text{id}_{\Cal V}+\psi(L_i)\ot\rho(\omega)$

\item{}
\qquad of $\W^\vn(1,\bold 1)$ over vector space $\frak V\ot_{\Cal K} \Cal V$.

\item{}
$(\text{\rm iii})$ \ Assume that the represetations $\phi$ and $\rho$ with a fixed point
$\omega\in$ $\Cal H^\vn(1)$ are given,
$\{\psi_a\}$ are

\item{}
\qquad a family of
representations of $\W^\vn(1,\bold 1)$ with the same representation vector space
$\frak V$ as $\phi$

\item{}
\qquad and satisfy $(31)$ relevant to $\phi$,
then the triples $(\phi,\;\sum_{a}k_a\psi_a;\;\rho)$
provide a family of repre-

\item{}
\qquad sentations of \ $\W^\vn(1,\bold 1)$ over
\ $\frak V\ot_{\Cal K} \Cal V$, which forms a vector space structure relative
to a

\item{}
\qquad fixed pair $(\phi,\rho)$.
\endproclaim
\demo{Proof} (i) Note that $\vn^{i+1}\psi(L_i)\psi(L_j)-
\vn^{j+1}\psi(L_j)\psi(L_i)=0$, it is clear that the pairs
$(\phi_{a\psi},\psi)$ still satisfy (31). That $\phi_{a\psi}$ is still
a representation of $\W^\vn(1,\bold 1)$ over $\frak V$ follows similarly
from the proof of next step (ii).

(ii) We need to show
$$
\vn^{i+1}\tilde \phi(e_{(i)})\tilde \phi(e_{(j)})-
\vn^{j+1}\tilde \phi(e_{(j)})\tilde \phi(e_{(i)})
=\tilde \phi(\{e_{(i)},e_{(j)}\}_\vn).
$$

\noindent
Note that
$\vn^{i+1}\binom{i+j}{i+1}_\vn-\vn^{j+1}\binom{i+j}{j+1}_\vn
=\binom{i+j+1}{i+1}_\vn-\binom{i+j+1}{j+1}_\vn$,
from (31), (13) and
$$\split
&\tilde \phi(e_{(i)})\tilde \phi(e_{(j)})\\
&=\bigl(\phi(e_{(i)})\ot\text{id}_{\Cal V}+\psi(L_i)\ot\rho(\omega)\bigr)\cdot\\
&\quad\ \bigl(\phi(e_{(j)})\ot\text{id}_{\Cal V}+\psi(L_j)\ot\rho(\omega)\bigr)\\
&=\phi(e_{(i)})\phi(e_{(j)})\ot\text{id}_{\Cal V}
+\phi(e_{(i)})\psi(L_j)\ot\rho(\omega)\\
&\quad +\psi(L_i)\phi(e_{(j)})\ot\rho(\omega)
+\psi(L_i)\psi(L_j)\ot\rho(\omega)^2,
\endsplit
$$
it is easy to obtain the required result:
$$\split
&\vn^{i+1}\tilde \phi(e_{(i)})\tilde \phi(e_{(j)})
-\vn^{j+1}\tilde \phi(e_{(j)})\tilde \phi(e_{(i)})\\
&=\left[\binom{i{+}j{+}1}{i{+}1}_\vn{-}\binom{i{+}j{+}1}{j{+}1}_\vn\right]
\bigl(\phi(e_{(i{+}j)})\ot\text{id}_{\Cal V}
{+}\psi(L_{i{+}j})\ot\rho(\omega)\bigr)\\
&=\left[\binom{i{+}j{+}1}{i{+}1}_\vn{-}\binom{i{+}j{+}1}{j{+}1}_\vn\right]
\tilde \phi(e_{(i+j)})\\
&=\tilde \phi(\{e_{(i)},e_{(j)}\}_\vn).
\endsplit
$$

(iii) If $\psi_a$ and $\psi_b$ satisfy (31) relevant to $\phi$,
then the pairs $(k\,\psi_a,\phi)$ $(k\in\Cal K$) and $(\psi_a+\psi_b,\phi)$ also
satisfy the same relation (31). So (iii) holds, thanks to (ii).

Thus we complete the proof.
\hfill\qed
\enddemo

\proclaim{Corollary 6.2}
A similar conclusion holds for the $q$-Witt algebra $\W^q(1)$
in the generic case, if $\W^\vn(1,\bold 1)$,
$\,_1\Cal K_\vn^l$, $\Cal H^\vn(1)$ and $\vn$ in $(31)$ are replaced by
$\W^q(1)$, $\,_1\Cal K_q^\infty$, $\Cal H^q$ and $q$, respectively.  \hfill\qed
\endproclaim

\vskip0.25cm
Here are two examples
to show some relationships between the universal truncated coinduced modules
and those given via the $\vn$-holomorph structure in Theorem 6.1.

\vskip0.2cm
\noindent
{\it Example 6.3} \ For $\W^\vn(1,\bold 1)$, if we take
$\frak V=\frak A^\vn(1,\bold 1)$, $\omega=e_{(0)}$, $\Cal V=\Cal K$ with
$\rho(e_{(0)})=t\in\Cal K$, and $\psi(L_i)=x^{(i)}\tau_\vn$, then we get the
realization given in section 5. It is a graded representation of
$\W^\vn(1,\bold 1)$ with a base space $\Cal K(t)=\Cal K$ over
$\frak A^\vn(1,\bold 1)\ot_{\Cal K}\Cal K(t)$, where
$\bigl(\frak A^\vn(1,\bold 1)\ot_{\Cal K}\Cal K(t)\bigr)_{(0)}\cong\Cal K(t)$.
Furthermore, if we take $\psi(L_i)=k\,x^{(i)}\tau_\vn, \forall k\in\Cal K$,
then we get all graded $\W^\vn(1,\bold 1)$-representations. In fact,
the universal truncated coinduced $\W^\vn(1,\bold 1)$-modules
only yield the special case where the enlarging coefficient $k=1$.
In the modular case, the same fact exists (see \cite{11, 13}).

\vskip0.2cm
\noindent
{\it Example 6.4} \
For the generic case in Corollary 6.2, if we take $\frak V=\frak A^q(1)$,
$\Cal V=\Cal K$, $\psi(L_i)=k\,x^{(i)}\tau_q$, $\omega=e_{(0)}$ and let $q$
tend to $1$, then we obtain every graded representation of the usual
(infinite dimensional) Witt algebra $\W(1)$ of Cartan type (cf. \cite{11, 13}).

\vskip0.2cm
\noindent
{\it Remark.}
Actually if we choose an arbitrary $\omega(\ne e_{(0)})\in \W^\vn(1,\bold 1)$
(resp. $\W^q(1)$), we arrive at some non-graded representations.
From Theorem 6.1 and Corollary 6.2, the $q$-holomorph structure
does lead to some new knowledge on representations of $\W^\vn(1,\bold 1)$
or $\W^q(1)$.

\vskip0.5cm
\noindent
{\it Acknowledgements.} \, I would like to express my sincere gratitude to Professor H. Strade for
the hospitality and support during my stay at Hamburg. I would thank
Professors Guangyu Shen and H. Strade for their helpful comments and valuable
suggestions. I also acknowledge the Alexander von Humboldt Stiftung for the support of
a Humboldt Research Fellowship.

\bye